\documentclass{article}
\usepackage[a4paper, left=0.4in, right=0.4in]{geometry}
\usepackage{authblk}

\usepackage{amsmath,amssymb,amsthm}
\usepackage{graphicx}
\usepackage{hyperref}

\usepackage{makecell}
\usepackage{lscape}
\usepackage[usenames]{color}
\definecolor{darkgreen}{rgb}{0,0.5,0}
\usepackage{amsfonts}

\usepackage{mathrsfs}
\usepackage{mathtools}
\usepackage{tikz}

\usepackage{subcaption}
\usepackage{multirow}

\newtheorem{lemma}{Lemma}[section]

\newtheorem{theorem}{Theorem}[section]
\newtheorem{definition}{Definition}[section]
\newtheorem{proposition}{Proposition}[section]
\newtheorem{example}{Example}

\title{A Fuzzy Geometric Study of Equidistant Sets in Fuzzy Metric Space}

\author[1]{Biswajit Singha}
\author[2]{Ronald Manríquez}
\author[2]{Cristian Carvajal}
\author[1]{Debjani Chakraborty}

\affil[1]{Department of Mathematics, Indian Institute of Technology Kharagpur, West Bengal, India}
\affil[2]{Laboratorio de investigaci\'on Lab[e]saM, Departamento de Matem\'atica,
 F\'isica y Computaci\'on, Universidad de Playa Ancha, Valpara\'iso, Chile.}

\begin{document}
\maketitle

\begin{abstract}

In this paper, the fuzzy Hausdorff distance is studied, and also the fuzzy equidistant set for two points of a fuzzy metric space is introduced. Here, the fuzzy metric space has been redefined using recently developed fuzzy geometry, and the equidistant sets have been constructed for two different fuzzy points. Different cases for the equidistant sets have been studied, considering two fuzzy points with separate spreads, externally tangent spreads, partially overlapping spreads, internally tangent spreads, fully overlapping spreads, and sets that coincide with the cores of fuzzy points.  The proposed construction provides a graded equidistant set that aligns with the classical midset when the metric is precise. Suitable numerical and pictorial examples are given to support the discussions and studies.
\end{abstract}


\section{Introduction}\label{sec_intro}
This paper aims to introduce the concept of fuzzy equidistant sets on a fuzzy metric space. 
 The concept of fuzzy metric space has been developed from multiple viewpoints. For instance, the George-Veeramani scheme, where a map $M: X\times X\times(0,\infty)\to(0,1]$ with a continuous $t$-norm encodes scale-dependent nearness \cite{georgeResultsFuzzyMetric1994}; the Kramosil-Michálek view, treating distance as a fuzzy subset of $\mathbb{R}_0^+$ \cite{kramosilFuzzyMetricsStatistical1975}; and their probabilistic counterparts, where distribution functions capture proximity \cite{schweizerStatisticalMetricSpaces1960}. Kaleva-Seikkala gives a more natural definition of fuzzy metric space \cite{kalevaFuzzyMetricSpaces1984}, where uncertainty is due to fuzziness rather than randomness. Further variants include intuitionistic and $L$-fuzzy metrics (\cite{parkIntuitionisticFuzzyMetric2004} and \cite{saadatiLfuzzyTopologicalSpaces2008}), and extensive literature, highlighting the works  \cite{dengFuzzyPseudometricSpaces1982}, \cite{puriDifferentialsFuzzyFunctions1983}, and \cite{shiGeneralizationsFuzzyMetric2023}.

In this paper, a fuzzy metric space in the sense of George-Veeramani \cite{georgeResultsFuzzyMetric1994} is considered due to its well-behaved (Hausdorff, metrizable) topology. The developed notions of fuzzy geometry \cite{chakraborty2024introduction},\cite{chakraborty2014analytical},\cite{ghosh2012analytical}, and \cite{ghosh2016analytical} are used to establish the closed form of proposed ideas. Fuzzy distance between two fuzzy points defined by Ghosh and Chakraborty \cite{ghosh2019introduction} is considered for defining fuzzy metric space in the sense of Kaleva-Seikkala \cite{kalevaFuzzyMetricSpaces1984}. Moreover, by providing canonical embeddings for crisp metrics, they allow geometric methods to be applied directly, while remaining more analytically tractable than other available approaches.

On the other hand, in a (crisp) metric space $(X,d)$, the equidistant set (or midset) of nonempty $A, B\subset X$ is defined to be $E(A,B)=\{x\in X:\, d(x,A)=d(x,B)\},$ where $d(x,A)=\inf~\{\textup{dist}(x,a):a\in A\}$. The sets $A$ and $B$ are called focal sets. $E(A,B)$ has been studied widely in \cite{lovelandWhenMidsetsAre1976} and \cite{wilkerEquidistantSetsTheir1975}. Much work has been done on this set because it generalizes the idea of a bisector or mediatrix. See, for example, \cite{foxEquidistantSetsAlexandrov2023},  \cite{ponceReflectionsEquidistantSets2022}, \cite{ponceEquidistantSetsGeneralized2014}, and \cite{zlepaloEquidistantSetsConic2019}.

In the fuzzy context,  \cite{manriquez2025approach} introduced the idea of fuzzy midset as a collection of fuzzy points, but the metric used does not allow for the study of topological characteristics of this set. Here, it was claimed that the cores of the fuzzy points belong to the midset defined by the cores of the fuzzy focal sets. However, it does not clarify where the remaining $\alpha$-cuts are situated. Also, it has not been discussed properly how the equidistant set behaves or is formed when the fuzzy focal sets are separate, partially overlapping, overlapping, or coincident with unequal spreads.
In the present study, this gap has been addressed. Suitable numerical and pictorial examples support the discussions and studies. The concept of a fuzzy metric space has been reintroduced, and the fuzzy equidistant set has also been discussed. Throughout the study, it is assumed that the compactness of the fuzzy focal sets and the convexity of fuzzy sets are satisfied.

The paper is structured as follows: Section~\ref{sec_preli} presents the fundamental prerequisites essential for the entirety of the study.  Section~\ref{sec_FMS_FHD} proposes the Fuzzy Hausdorff distance and the fuzzy metric space. Section~\ref{sec_Redefine_proposal} introduces the central concept of this work, the fuzzy equidistant set. In section~\ref{sec_Redefine_proposal}, it has been shown that with the proposed idea of fuzzy metric space, fuzzy equidistant set is invariant for both the metric spaces $\bigl(X,\widetilde{M}_d,*\bigl)$ and $\bigl(X,\widetilde{d}, L, R\bigr)$. A natural fuzzy generalization of $\widetilde{E}(\widetilde{A},\widetilde{B})$ has also been proposed by establishing it with a membership function. Finally, in section~\ref{sec_conclusion}, some comments and conclusions are presented.

\section{Preliminaries}\label{sec_preli}
In this section, a summary of definitions and results has been given.  Mainly, notions about fuzzy distance and fuzzy metric spaces in the sense of George and Veeramani \cite{georgeResultsFuzzyMetric1994} are given. The fundamental components of fuzzy geometry, as outlined in \cite{chakraborty2014analytical}, have also been summarised.

\begin{definition}[Fuzzy set, \cite{zadehFuzzySets1965}]
Let $X$ be a classical set. Then, the set of ordered pairs 
$$\tilde{A}=\{\left(x,\mu_{\tilde{A}}(x)\right) : x \in X, \,\ \mu_{\tilde{A}}(x) \in [0,1]\}$$
is called a fuzzy set on $X$. The evaluation function $\mu_{\tilde{A}}(x)$ is called the membership function or the grade of membership of $x$ in $\tilde{A}$.
\end{definition}

\begin{definition}[$\alpha$-cut, \cite{duboisFuzzySetsSystems1980}]
For a fuzzy set $\tilde{A}$, its $\alpha$-cut is denoted by $\tilde{A}^\alpha$ and is defined by
$$\tilde{A}^{\alpha} =\begin{cases} 
\{x\in X : \; \mu_{\tilde{A}}(x) \geq \alpha\}&\text{if}\quad 0<\alpha\leq1\\
\textup{cl}\{x\in X : \; \mu_{\tilde{A}}(x) >\alpha\}&\text{if}\quad \alpha=0,
\end{cases}$$ 
where $\textup{cl}\{x\in X : \; \mu_{\tilde{A}}(x) >\alpha\}$ denotes the topological closure.
\end{definition}

\begin{definition}[Fuzzy numbers, \cite{buckley1997fuzzy}] 
    A fuzzy set $\widetilde{A}$ of $\mathbb{R}$ is called a fuzzy number if its membership function $\mu$ has the following properties:

    \begin{enumerate}
        \item[(i)] $\mu(x|\widetilde{A})$ is upper semi-continuous,
        \item[(ii)]$\mu(x|\widetilde{A})=0$ outside some interval $[a,d]$, and
        \item [(iii)]
        there exist real numbers $b$ and $c$ so that $a\leq b\leq c\leq d$ and $\mu(x|\widetilde{A})$ is increasing on $[a,b]$ and decreasing on $[c,d]$, and $\mu(x|\widetilde{A})=1$ for each $x \in [b,c]$.
    \end{enumerate}    
\end{definition}
If $b=c$ and let $f(x)=\mu(x|\widetilde{A}) ~~\forall x \in [a,b]$ and $g(x)=\mu(x|\widetilde{A}) ~~\forall x \in [b,d]$ are linear functions, then the fuzzy number is called triangular fuzzy number and is denoted by $(a,b,d)$ or $(a/b/d)$.

\begin{definition}[Fuzzy number along a line,\cite{chakraborty2014analytical} \cite{ghosh2012analytical} \cite{ghosh2019introduction}]
In defining a fuzzy number, conventionally a real line $(\mathbb{R})$ is taken as the universal set. Accordingly a fuzzy number on the real line, $\widetilde{p}$, say, is defined as a collection of real numbers
those are neighboring to the real number p with different membership values. Instead of a real line, consider any line on the plane $\mathbb{R}^2$ then $\widetilde{p}$ can be represented as a fuzzy number on that line. On the $x$-axis, the membership function of $\widetilde{p}$ may be written as $\mu\bigl((x,0)| ~\widetilde{p})\bigr)=\mu(x|\widetilde{p})~~ \forall x \in \mathbb{R}$. More explicitly:   
$\mu\bigl((x,0)|\widetilde{p})\bigr) =\begin{cases} 
\mu(x|\widetilde{p})&\text{if}\quad y=0\\
0&\text{elsewhere}.
\end{cases}$

Instead of $x$-axis, consider any line $ax+by=c$ on the $xy$-plane as the universal set, which makes an angle $\theta$, then a bijective transformation $T:\mathbb{R}^2 \to \mathbb{R}$ is defined which is of the form:
\begin{equation}
\label{transformation}
T(x,y)=\Bigl((x-h) \cos{\theta} + (y-k) \sin{\theta},~-(x-h) \sin{\theta}+(y-k) \cos{\theta}\Bigr)
\end{equation}
transforms $x$-axis to the line $ax+by=c$, then
\begin{itemize}
\item[(i)] the point $(h,k)=\Bigl(\frac{ac}{a^2+b^2},~\frac{bc}{a^2+b^2}\Bigr)$ being point of intersection for $ax+by~=c$ and its perpendicular line through origin $bx-ay=c$, and
\item[(ii)] $\theta$ being the angle of elevation of $ax+by~=c$ with the positive $x$-axis.
\end{itemize}
So if $\widetilde{p}$ is a fuzzy number on the $x$-axis, then on the line $l_{p\theta}$ (say) $: ax+by=c$ the fuzzy number can be considered as $T(\widetilde{p})$.
\end{definition}
\begin{definition}
[Fuzzy Point, \cite{ghosh2019introduction}]
    Let $(a,b) \in \mathbb{R}^{2}$, a fuzzy set in $\mathbb{R}^{2}$ is called a fuzzy point $\widetilde{P}(a,b)$ if its membership function satisfies the following properties:
    \begin{itemize}
        \item[(i)] $\mu((x,y)| \widetilde{P}(a,b))$ is upper semi-continuous,
        \item[(ii)] $\mu((x,y)| \widetilde{P}(a,b)) = 1$ only at (x,y) = (a,b), and
        \item[(iii)] $\widetilde{P}^\alpha(a,b)$ (alpha-cut) is a convex and compact subset of $\mathbb{R}^{2}$ for all $\alpha \in [0,1].$
    \end{itemize}
    \label{fuzzyPoint}
\end{definition}

\begin{definition}
    [Same points, \cite{ghosh2019introduction}]
    Two points $(x_{1}, y_{1})$, $(x_{2}, y_{2})$ on the supports of the fuzzy points $\widetilde{P}$(a,b), $\widetilde{P}$(c,d) respectively are said to be same points if
    \begin{itemize}
        \item[(i)]  $\mu((x_{1}, y_{1})| \widetilde{P}(a,b))$ = $\mu((x_{2}, y_{2})| \widetilde{P}(c,d))$, 

        \item[(ii)] $(x_{1}, y_{1})$ and $(x_{2}, y_{2})$ lie on the same side of the line joining $(a,b)$ and $(c,d)$, and

        \item[(iii)] the line joining $(x_{1}, y_{1})$ to $(a,b)$ and the line joining $(x_{2}, y_{2})$ to $(c,d)$ are parallel.
    \end{itemize}
\label{same}
\end{definition}
\begin{definition}
    [Inverse points, \cite{ghosh2019introduction}] 
    Two points $(x_{1}, y_{1})$, $(x_{2}, y_{2})$ on the supports of the fuzzy points $\widetilde{P}$(a,b), $\widetilde{P}$(c,d) respectively are said to be inverse points if
    \begin{itemize}
        \item[(i)]  $\mu((x_{1}, y_{1})| \widetilde{P}(a,b))$ = $\mu((x_{2}, y_{2})| \widetilde{P}(c,d))$, 
        \item[(ii)] $(x_{1}, y_{1})$ and $(x_{2}, y_{2})$ lie on different sides of the line joining $(a,b)$ and $(c,d)$, and
        \item[(iii)] the line joining $(x_{1}, y_{1})$ to $(a,b)$ and the line joining $(x_{2}, y_{2})$ to $(c,d)$ are parallel.
    \end{itemize}
\label{inverse}
\end{definition}
\begin{definition}
   [Fuzzy distance between fuzzy points, 
\cite{chakraborty2014analytical} \cite{ghosh2012analytical} \cite{ghosh2019introduction}]
Let $\widetilde{A}$ and $\widetilde{B}$ be two fuzzy points, then the Fuzzy distance is redefined as follows
\begin{center}    $\mu(\widetilde{d^\alpha}~|~\widetilde{d})=sup\{\alpha: d=d(u,v),
     \forall $$(u,v)$ \text{lies on the line joining inverse points}   \\ $ u \in \widetilde{A}^0$  \text{and}  $v \in \widetilde{B}^0$  \text{and}  $\mu(u\mid \widetilde{A})=\mu(v\mid \widetilde{B})=\alpha, d \text{ is euclidean distance}\big\}$\\
     \text{and}\\ 
     $\widetilde{d}(\widetilde{A},\widetilde{B}) = \bigcup_{\alpha=0}^{1}  \widetilde{d}^{\alpha}\bigl(\widetilde{A}^\alpha,\widetilde{B}^\alpha\bigl)$
\end{center}
Here, $\widetilde{d}(,)$ is the fuzzy distance metric and $\widetilde{A}^0$, $\widetilde{B}^0$ are supports of fuzzy points $\widetilde{A}$, $\widetilde{B}$ respectively. 
\label{def_fuzzy_distance} 
\end{definition}
\begin{example}   
Let $\widetilde{A}(a_1,a_2)$ and $\widetilde{B}(b_1,b_2)$ be two fuzzy points with elliptical bases  having spreads $S_1 = (p_{1},p_{2})$ and $S_2 = (q_{1},q_{2})$. In every $\alpha$-cut ($\alpha \in [0,1]$)  consider a line $l_{\theta}$ which makes an angle $\theta$ with real axis. The intersecting points on $\widetilde{A}^{\alpha}$ are $\underline{A}^{\alpha}_{\theta}$ and $\overline{A}^{\alpha}_{\theta}$ respectively. Similarly, the intersecting points on $\widetilde{B}^{\alpha}$  are $\underline{B}^{\alpha}_{\theta}$ and $\overline{B}^{\alpha}_{\theta}$. Then,
\begin{equation}
\begin{aligned}
\underline{A}^{\alpha}_{\theta} &= \bigl(a_1 - p_1*(1-\alpha)*\cos(\theta),\, a_2 - p_2*(1-\alpha)*\sin(\theta)\bigr) \\
\overline{A}^{\alpha}_{\theta} &= \bigl(a_1 + p_1*(1-\alpha)*\cos(\theta),\, a_2 + p_2*(1-\alpha)*\sin(\theta)\bigr) \\
\underline{B}^{\alpha}_{\theta} &= \bigl(b_1 - q_1*(1-\alpha)*\cos(\theta),\, b_2 - q_2*(1-\alpha)*\sin(\theta)\bigr) \\
\overline{B}^{\alpha}_{\theta} &= \bigl(b_1 + q_1*(1-\alpha)*\cos(\theta),\, b_2 + q_2*(1-\alpha)*\sin(\theta)\bigr)
\end{aligned}
\label{eq_alphacuts}
\end{equation}
where $\theta \in [0,\pi]$. The distance $\widetilde{d}^{\alpha}$ is defined using definition \ref{def_fuzzy_distance}
\begin{equation}       \widetilde{d}^{\alpha}\bigl(\widetilde{A}^\alpha,\widetilde{B}^\alpha\bigl) = \biggl(\min_{\theta \in [0,\pi]} \underline{\lambda}, d((a_1,a_2),(b_1,b_2)),\max_{\theta \in [0,\pi]} \overline{\lambda}\biggl)
       \label{eq_fuzzy_dist}
\end{equation}
where,   $
\underline{\lambda} =
\min~\!\Bigl( d~\!\bigl(\underline{A}^{\alpha}_{\theta},\, \overline{B}^{\alpha}_{\theta}\bigr), \;
              d~\!\bigl(\overline{A}^{\alpha}_{\theta},\, \underline{B}^{\alpha}_{\theta}\bigr) \Bigr)
$   and   $
\overline{\lambda} = 
\max~\!\Bigl( d~\!\bigl(\underline{A}^{\alpha}_{\theta},\, \overline{B}^{\alpha}_{\theta}\bigr), \;
              d~\!\bigl(\overline{A}^{\alpha}_{\theta},\, \underline{B}^{\alpha}_{\theta}\bigr) \Bigr).
$ $\widetilde{A}^\alpha,~\widetilde{B}^\alpha$ are alpha-cuts  and $d(,)$ denotes the usual euclidean distance. Considering all $\alpha$-cuts the fuzzy distance between fuzzy points $\widetilde{A}(a_1,a_2)$ and $\widetilde{B}(b_1,b_2)$ is 
\begin{equation}
    \widetilde{d}(\widetilde{A},\widetilde{B}) = \bigcup_{\alpha=0}^{1}  \widetilde{d}^{\alpha}\bigl(\widetilde{A}^\alpha,\widetilde{B}^\alpha\bigl)
    \label{eq_dist}
\end{equation}
$\widetilde{d}(\widetilde{A},\widetilde{B})$ is a fuzzy set with membership function as follows:
\begin{equation}
\mu_{\widetilde{d}}(x) = \begin{cases}
L(d) = {\min}\{1 - \phi(x,\theta)\}, ~~~~ \text{for } \theta_{\max},\\
R(d) = {\max}\{1 -  \phi(x,\theta)\},~~~~\text{for } \theta_{\min},
\end{cases}
\end{equation}
where, 
\begin{equation}
\begin{aligned}
&\phi(x,\theta) = \frac{-\bigl(d_1 R_1 \cos\theta + d_2 R_2 \sin\theta\bigr)\pm \sqrt{\bigl(d_1 R_1 \cos\theta + d_2 R_2 \sin\theta\bigr)^2 - (R_1^2\cos^2\theta + R_2^2\sin^2\theta)(d_c^2 - x^2)}}{(R_1^2\cos^2\theta + R_2^2\sin^2\theta)}\\
& \text{and }R_1 = p_1 + q_1, ~~ R_2 = p_2 + q_2, ~~
d_1 = a_1 - b_1,~~ d_2 = a_2 - b_2,
~~
d_c = \sqrt{d_1^2 + d_2^2} 
\end{aligned}
\end{equation}
\label{elli_pt_fuz_dis}
\end{example}
\begin{example}
    Let $\widetilde{A}(1,0)$  and $\widetilde{B}(5,2)$ (Fig. \ref{fp_dist_eg2d}) be two fuzzy points with membership functions as 
    \begin{center}
        $\mu_{1}((x,y)|\widetilde{A}(1,0)) = \begin{cases}
        1-\sqrt{(x-1)^{2}+(y-0)^{2} } & if\  (x-1)^{2}+(y-0)^{2} \leq 1 \\
        0 & otherwise \\
    \end{cases}$
    \end{center}

    \begin{center}
        $\mu_{2}((x,y)|\widetilde{B}(5,2)) = \begin{cases}
        1-\sqrt{(\frac{x-5}{1})^{2}+(\frac{y-2}{1.5})^{2}} & if\  (\frac{x-5}{1})^{2}+(\frac{y-2}{1.5})^{2}\leq 1\\
        0 & otherwise \\
    \end{cases}$
    \end{center}

\begin{figure}[h]
\centering
\begin{subfigure}{0.45\textwidth}
\includegraphics[width=0.9\linewidth]{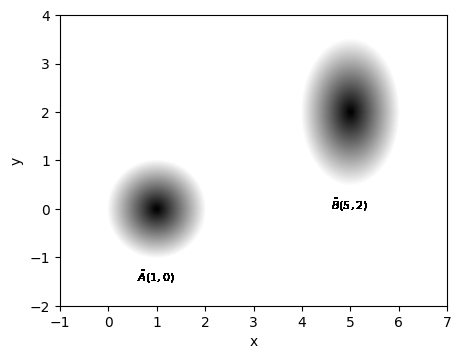} 
\caption*{\small \parbox{6cm}{Fuzzy points $\widetilde{A}(1,0)$ and $\widetilde{B}(5,2)$ with different spreads}}
\label{fuzzypoints}
\end{subfigure}
\begin{subfigure}{0.45\textwidth}
\includegraphics[width=0.9\linewidth]{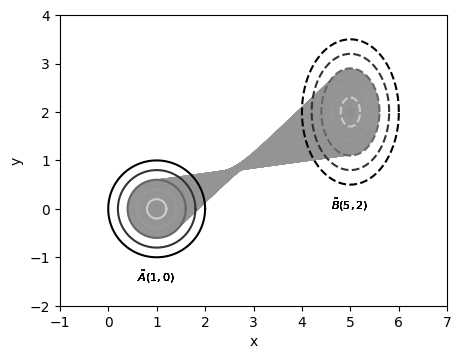}
\caption*{\small \parbox{5cm} {Fuzzy distance at $\alpha=0.4$-cut}}
\label{0.4_cut}
\end{subfigure}
\caption{Fuzzy points and fuzzy distance}
\label{fp_dist_eg2d}
\end{figure}
 Here, $\min_{\theta \in [0,\pi]} \underline{\lambda}$ occurs at $\theta = 0.4631$ and $\max_{\theta \in [0,\pi]} \overline{\lambda}$ is at $\theta = 3.8168$. Using Eq.\ref{eq_fuzzy_dist}, the distance at $\alpha$-cut is as follows: 
\begin{equation*}
\scalebox{0.92}{$   \widetilde{d}^{\alpha}\bigl(\widetilde{A}^\alpha,\widetilde{B}^\alpha\bigl) = \Bigl(\sqrt{5.667025+9.883959\alpha+4.449017\alpha^2},\\ 4.472136,\sqrt{43.618887-28.497955\alpha+4.879067\alpha^2}\Bigl)
$}
\end{equation*}
Considering all $\alpha$-cuts together, the fuzzy distance between fuzzy points $\widetilde{A}(1,0)$ and $\widetilde{B}(5,2)$ is\\
\begin{center}
 $\widetilde{d}(\widetilde{A},\widetilde{B}) = (2.380551, 4.472136,6.604459)$.  
\end{center}

\label{fuzzy_dist_ex}
\end{example}

Instead of considering two fuzzy points with different bases, if two fuzzy points have similar bases, then the following proposition applies.

\begin{proposition}
If $\widetilde{A}$ and $\widetilde{B}$ are two fuzzy points with circular bases, then the values $\underset{\theta \in [0,\pi]}{min}\underline{\lambda}$ and $\underset{\theta \in [0,\pi]}{max}\overline{\lambda}$ (see Eq.\ref{eq_fuzzy_dist}) are attained at the angle whose value coincides with the slope of the line segment connecting the cores of the fuzzy points $\widetilde{A}$ and $\widetilde{B}$.
\label{cor1}
\end{proposition}
\begin{figure}[ht]
    \centering   \includegraphics[width=0.55\linewidth]{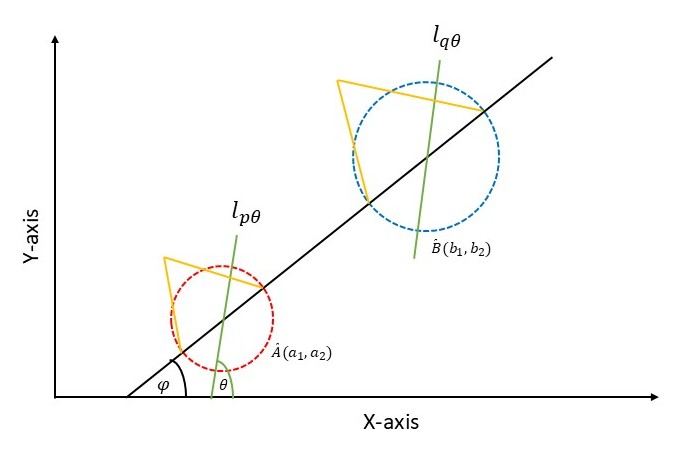}
    \caption{$l_{p\theta}$ and $l_{q\theta}$ make $\theta$ with $x$-axis and line joining $(a_1,a_2)$ and $(b_1,b_2)$ makes $\psi$ with $x$-axis}
    \label{fig:theta_psi}
\end{figure}

\textbf{Proof :} Let $\widetilde{A}(a_1,a_2)$ and $\widetilde{B}(b_1,b_2)$ be two fuzzy points with circular bases having spreads $S_1=r_1$ and $S_2=r_2$. Then using equation \ref{eq_alphacuts} we have,
\begin{equation}
\begin{aligned}
d^2\bigl(\overline{A}^{\alpha}_{\theta},\underline{B}^{\alpha}_{\theta}\bigr)
&= d_c^2 + (r_1+r_2)^2(1-\alpha)^2 
+ 2(r_1+r_2)(1-\alpha)\bigl[(a_1-b_1)\cos{\theta} + (a_2-b_2)\sin{\theta}\bigr] \\
d^2\bigl(\underline{A}^{\alpha}_{\theta},\overline{B}^{\alpha}_{\theta}\bigr)
&= d_c^2 + (r_1+r_2)^2(1-\alpha)^2 
- 2(r_1+r_2)(1-\alpha)\bigl[(a_1-b_1)\cos{\theta} + (a_2-b_2)\sin{\theta}\bigr]
\end{aligned}
\label{distsq}
\end{equation}

where, $d_c^2 = (a_1-b_1)^2+(a_2-b_2)^2$.\\

Differentiating Eq.\ref{distsq} with respect to $\theta$, we have $2(r_1+r_2)(1-\alpha)[(a_1-b_1)\sin{\theta}-(a_2-b_2)\cos\theta]=0$. For all the values of $\alpha \in [0,1)$ the minimum value of $\theta$ is $\tan^{-1}[\frac{a_2-b_2}{a_1-b_1}]$. So when $l_{p\theta}$ and $l_{q\theta}$ both are at the angle of $\theta_{min} = \tan^{-1}[\frac{a_2-b_2}{a_1-b_1}]$ we get the fuzzy distance between two fuzzy points. For $\alpha=1$, $\widetilde{A}$ and $\widetilde{B}$ becomes non fuzzy. So no need to calculate the $\alpha$-cuts. 
Again, the slop of the line $\frac{x-a_2}{a_1-a_2}=\frac{y-b_2}{b_1-b_2}$, joining $(a_1,a_2)$ and $(b_1,b_2)$ is  $\tan{\psi}=\frac{a_2-b_2}{a_1-b_1}$ that is $\psi= \tan^{-1}[\frac{a_2-b_2}{a_1-b_1}]$. Which is the same as $\theta_{min}$. 

This shows that inverse points for the fuzzy distance between two fuzzy points with circular bases occur at the angle of the line joining the cores of the fuzzy points. This guarantees us that for fuzzy distance between two fuzzy points is the distance between the fuzzy numbers along the direction which makes $\psi$. If we consider two fuzzy points with elliptical bases with spreads $S_1=(p_1,p_2)$ and $S_2=(p_2,p_1)$, then we can obtain the same result.

\begin{definition}[Fuzzy distance between a fuzzy point and a fuzzy set, \cite{manriquez2025approach}]
Let $\widetilde{A}$ be a fuzzy set and $\widetilde{P}$ be a fuzzy point then fuzzy distance denoted by $\widetilde{d}$, between $\widetilde{P}$ and $\widetilde{A}$ is defined as $\widetilde{d}(\widetilde{P},\widetilde{A})=\widetilde{d}(\widetilde{P}_A,\widetilde{a})$
\label{fuzzydist_set_pt}
\end{definition}

\begin{definition}[Fuzzy midset, \cite{manriquez2025approach}]
Let $\widetilde{A}$ and $\widetilde{P}$ be a fuzzy sets on $\mathbb{R}^2$. The fuzzy midset between $\widetilde{A}$ and $\widetilde{B}$, denoted by $\widetilde{\mathcal{M}}(\widetilde{A},\widetilde{B})$ may be defined by the collection of fuzzy points $\widetilde{P}$, such that $\widetilde{d}(\widetilde{P},\widetilde{A})=\widetilde{d}(\widetilde{P},\widetilde{B})$ and $\widetilde{P}^1 \in \mathcal{M}\bigl(\widetilde{A}^1,\widetilde{B}^1\bigr)$.
\label{f_midset}
\end{definition}
\begin{definition}[$t$-norm, \cite{schweizerStatisticalMetricSpaces1960}]
A binary operation $*:[0,1]\times[0,1]\rightarrow[0,1]$ is called a continuous $t$-norm if $([0,1],*)$ is an Abelian topological monoid with unit 1, such that $a*b\leq c*d$ whenever $a\leq b$ and $c\leq d$, $a,b,c,d\in[0,1]$.
\end{definition}
\begin{definition}[Fuzzy metric space, \cite{georgeResultsFuzzyMetric1994}\cite{kramosilFuzzyMetricsStatistical1975}]\label{def_fmetric}

A fuzzy metric space is an ordered triple $(X,M,*)$ such that $X$ is a nonempty set, $*$ is a continuous $t$-norm and $M$ is a fuzzy set of $X^2\times (0,\infty)$ satisfying the following conditions, for all $x,y,z\in X$, $s,t>0$:
\begin{enumerate}
\item $M(x,y,t)>0$,
\item $M(x,y,t)=1$ if and only if $x=y$,
\item $M(x,y,t)=M(y,x,t)$;
\item $M(x,y,t)*M(y,z,s)\leq M(x,z,t+s)$,
\item $M(x,y,\cdot):(0,\infty)\rightarrow [0,1]$ is continuous.
\end{enumerate}
If  $(X,M,*)$ is a fuzzy metric space, we say that $(M,*)$ is a fuzzy metric on $X$.
\end{definition}
Let $x\cdot y$ be the usual multiplication for all $x,y\in[0,1]$, and let $M$ be the function defined on $X^2\times (0,\infty)$ such that
\begin{equation} \label{eq:fuzzy_metric}
M(x,y,t) = \frac{t}{t + d(x,y)}, \quad x,y \in X,\ t > 0.
\end{equation}
 Then, $\left(X,M,\cdot\right)$ is a fuzzy metric space called standard fuzzy metric space (see \cite{georgeResultsFuzzyMetric1994}), and $\left(M,\cdot\right)$ will be called the standard fuzzy metric of $d$. In \cite{georgeResultsFuzzyMetric1994}, it was demonstrated that every fuzzy metric $(M,*)$ on $X$ generates a topology $\tau$ such that $(X,\tau)$ is a metrizable topological space. In metric space instead to distance $d(x,y)$, $M(x,y,t)$ is considered, which is the degree of similarity or closeness between \(x\) and \(y\) at time/distance \(t\) (often \(t\in [0,\infty )\) or \((0,\infty )\)).
 
 The following results are useful in the following sections.
 
\begin{proposition}[Proposition 1 in \cite{rodriguez-lopezHausdorffFuzzyMetric2004}]\label{prop_Mconti}]
Let $(X,M,*)$ be a fuzzy metric space. Then $M$ is a continuous function on $X^2\times (0,\infty)$.
\end{proposition}

\begin{definition}[Definition 2.4 in \cite{veeramaniBestApproximationFuzzy2001}]\label{def_xfromt}
Let $A$ be a nonempty subset of a fuzzy metric space  $(X,M,*)$. For $x\in X$ and $t>0$, let $M(x,A,t):=\sup_{a\in A}\{M(x,a,t)\}$.
\end{definition}

\begin{lemma}[Lemma 2 in \cite{rodriguez-lopezHausdorffFuzzyMetric2004}]
Let $(X,M,*)$ be a fuzzy metric space. Then, for each $a\in X$ and $B\in\mathcal{K}_0(X)$ the function $t\longmapsto M(a,B,t)$ is continuous on $(0,\infty)$, where $\mathcal{K}_0(X)$ is the set of nonempty compact subsets of $X$.
\end{lemma}

\begin{definition}[Fuzzy metric space \cite{kalevaFuzzyMetricSpaces1984}]\label{gen_FMS}

Let $X$ be a non-empty set, $d$ a mapping from $X \times X$ into $G$ (set of all 
non-negative fuzzy numbers of E) and let the mappings $L,R:[0,1] \times [0,1] \to [0,1]$ be symmetric, non-decreasing in both arguments and satisfy $L(0,0)=0$ and $R(1,1)=1$. Denote
\[
[d(x,y)]_{\alpha} = [\lambda_{\alpha}(x,y),\rho_{\alpha}(x,y)] ~\text{for}~ x,y \in X, ~0<\alpha\leq 1.
\]
The quadruple $(X,d,L,R)$ is called a fuzzy metric space and $d$ is a fuzzy metric, if
\begin{enumerate}
    \item $d(x,y) = \overline{0}$ if and only if $x=y$,
    \item $d(x,y) = d(y,x)$ for all $x,y \in X$,
    \item for all $x,y,z \in X$
    \begin{enumerate}
        \item $d(x,y)(s+t) \leq R(d(x,y)(s), d(z,y)(t))$
        \\
        whenever $s \leq \lambda_1(x,z)$, $t \leq \lambda_1(z,y)$ and $s+t \leq \lambda_1(x,y)$
        \item $d(x,y)(s+t) \geq L(d(x,y)(s), d(z,y)(t))$
        \\
        whenever $s \geq \lambda_1(x,z)$, $t \geq \lambda_1(z,y)$ and $s+t \geq \lambda_1(x,y)$
    \end{enumerate} \label{tri_ineq_Gen_FMS}
\end{enumerate}
\end{definition}
 \begin{theorem}
     In a fuzzy metric space $(X,d,L,R )$ the triangle inequality \ref{tri_ineq_Gen_FMS} in Definition \ref{gen_FMS} is equivalent to $d(x,y) \preceq d(x,y) + d(z,y)$ considering L as Min and R as Max. [see proof in \cite{kalevaFuzzyMetricSpaces1984}]
 \end{theorem}

\section{Introduction to Fuzzy Metric Space}\label{sec_FMS_FHD}
In this section,  the definition of the fuzzy metric space [see Definition \ref{def_fmetric}] incorporating fuzzy points has been proposed. Subsequently, the fuzzy distance between two fuzzy points [see Definition \ref{def_fuzzy_distance}] is proven to be a metric, and the definition of the fuzzy Hausdorff distance is also proposed.

\begin{definition}[Redefined Fuzzy Metric Space] \label{Redefined Fuzzy metric space}

A fuzzy metric space is an ordered triple $(X,\widetilde{M_d},*)$ such that $X$ is the universal set, $*$ is a continuous $t$-norm and $\widetilde{M}$ is a fuzzy set on $X^2\times (0,\infty)$ satisfying the following conditions, for all fuzzy points $\widetilde{A},\widetilde{B},\widetilde{C}\in X$, $t,s>0$:
\begin{enumerate}
\item $\widetilde{M_d}(\widetilde{A},\widetilde{B},t)\succ
0$,
\item $\widetilde{M_d}(\widetilde{A},\widetilde{B},t)\cong 1$ if and only if $\widetilde{A} \cong \widetilde{B}$,
\item $\widetilde{M_d}(\widetilde{A},\widetilde{B},t)\cong \widetilde{M_d}(\widetilde{B},\widetilde{A},t)$; 
\item $\widetilde{M_d}(\widetilde{A},\widetilde{B},t)*\widetilde{M}(\widetilde{B},\widetilde{C},s) \preceq \widetilde{M}(\widetilde{A},\widetilde{C},t+s)$,
\item $\widetilde{M_d}(\widetilde{A},\widetilde{B},\cdot):(0,\infty)\rightarrow [0,1]$ is continuous.
\end{enumerate}
Here, $\cong$ means 'almost equals' in fuzzy sense and $\preceq$ means 'fuzzily less equal to'. Moreover, $(X,\widetilde{M_d},*)$ is a fuzzy metric space, if $(\widetilde{M_d},*)$ is a fuzzy metric on $X$ and $\widetilde{d}(,)$ is the fuzzy distance (definition \ref{def_fuzzy_distance}).
\end{definition} 
All the conditions mentioned in the preceding Definition \ref{Redefined Fuzzy metric space} appear to have a remarkably similar interpretation to \cite{kramosilFuzzyMetricsStatistical1975}. These conditions also convey the notion that the characteristics of identity, non-negativity, and symmetry are subject to fuzziness or uncertainty. It is also to be noted that fuzzy metric spaces have fuzzy points as their elements, i.e., they are sets of fuzzy points \cite{xia2004fuzzy}.

For instance, a fuzzy metric $\widetilde{M_d}$ can be defined as follows. Let $\widetilde{M_d}$ be a function of two fuzzy points $\widetilde{A}$ and $\widetilde{B}$ on $X^2 \times (0, \infty)$ and can be established as follows:
\begin{equation} \label{eq:Re_def_fuzzy_metric}
\widetilde{M_d}(\widetilde{A},\widetilde{B},t) = \frac{t}{t + \widetilde{d}(\widetilde{A},\widetilde{B})}, \quad \widetilde{A},\widetilde{B} \in X,\ t>0 
\end{equation}
Here, $ \widetilde{M_d}(\widetilde{A},\widetilde{B},t)$ gives the degree to which the sets \( \widetilde{A} \) and \( \widetilde{B} \) are close to each other under the fuzzy metric \( \widetilde{M_d} \) at \( t \).

\begin{example}
    If we consider the example
    \ref{fuzzy_dist_ex}. For every $\widetilde{A}, \widetilde{B} \in X, t>0$
    \[
    \begin{aligned}    \widetilde{M}_d(\widetilde{A},\widetilde{B},t)
    &= \frac{t}{t + \widetilde{d}(\widetilde{A},\widetilde{B})}  \\
    &= \frac{t}{t + (2.380551, 4.472136,6.604459)} \\
    &= \frac{t}{ (t+2.380551, t+4.472136,t+6.604459)} \\
&= \Biggl(\frac{t(t+2.339813)}{(t+4.472136)^2},\frac{t}{t+4.472136},\frac{t(t+6.563721)}{t+4.472136}\Biggr)
    \end{aligned}
    \]

\begin{figure}[h]
    \centering    \includegraphics[scale=.35]{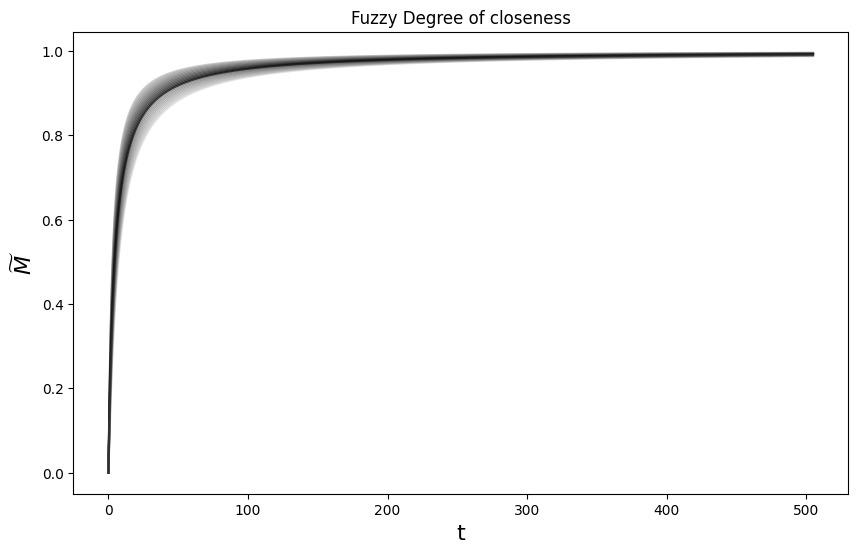}
    \caption{degree of closeness between the fuzzy points $\widetilde{A}(1,0)$, $\widetilde{B}(5,2)$ }
    \label{deg_colseness}
    \end{figure}
\label{ex_deg_closeness}
\end{example}
In this example, it is evident that the fuzziness or uncertainty is intrinsically linked to the degree of closeness. From Figure-\ref{deg_colseness}, it is observed that when $t$ is close to zero, the uncertainty is low. As $t$ increases, the uncertainty grows. For very large values of $t$, however, the uncertainty decreases again and eventually vanishes.  
Moreover, when $t$ is close to zero, $\widetilde{M_d}$ is also close to zero, and as $t$ becomes large, $\widetilde{M_d}$ approaches the value 1.

Normally, fuzzy distance shows us how far the fuzzy sets or fuzzy points $\widetilde{A}$ and $\widetilde{B}$ are apart. Thus, a smaller distance implies higher degree of closeness, while a larger distance means lower degree of closeness.\\

\begin{proposition}
    Fuzzy distance between two fuzzy points is fuzzy metric and  $(X,\widetilde{M_d},*)$ is a metric space. 
\end{proposition}
\textbf{Proof : }
From the equation \ref{eq_dist} we have $\widetilde{d}(\widetilde{A},\widetilde{B})=\bigcup_{\alpha=0}^{1}  \widetilde{d}^{\alpha}\bigl(\widetilde{A}^\alpha,\widetilde{B}^\alpha\bigl)=(\underline{\lambda},d(A,B),\overline{\lambda})$
\begin{enumerate}
    \item $\widetilde{d}(\widetilde{A},\widetilde{B}) \approx \widetilde{0}$ if and only if $\widetilde{A} \cong \widetilde{B}$\\
 Let, $\widetilde{A}\cong\widetilde{B}$ then $\widetilde{d}(\widetilde{A},\widetilde{A})= \bigcup_{\alpha=0}^1\widetilde{d}^{\alpha}\bigl(\widetilde{A}^{\alpha},\widetilde{A}^{\alpha}\bigl)=(\underline{\lambda},0,\overline{\lambda})\cong \widetilde{0}$

Conversely, let $\widetilde{d}(\widetilde{A},\widetilde{B})= 0$. Then, $\bigcup_{\alpha=0}^1\widetilde{d}^{\alpha}\bigl(\widetilde{A}^{\alpha},\widetilde{B}^{\alpha}\bigl) ~=~ 0,$ implies 
$\widetilde{d}^{\alpha}\bigl(\widetilde{A}^{\alpha},\widetilde{B}^{\alpha}\bigl)=0  \text{ for each }  \alpha \in [0,1].$

Thus, $\d(\underline{A}_\theta^\alpha,\underline{B}_\theta^\alpha)=0= d(\overline{A}_\theta^\alpha,\overline{B}_\theta^\alpha)$ and this implies that $\widetilde{A} \cong\widetilde{B}.$

\item $\widetilde{d}(\widetilde{A},\widetilde{B}) = \widetilde{d}(\widetilde{B},\widetilde{A})$ for all $\widetilde{A}, \widetilde{B} \in X$

 Proof is obvious as $d(\underline{A}_\theta^\alpha,\underline{B}_\theta^\alpha)$ and $d(\overline{A}_\theta^\alpha,\overline{B}_\theta^\alpha)$ are the euclidean distances, and thus symmetric.

\item $\widetilde{d}(\widetilde{A},\widetilde{B}) \preceq \widetilde{d}(\widetilde{A},\widetilde{C}) + \widetilde{d}(\widetilde{C},\widetilde{B})$ for all $\widetilde{A}, \widetilde{B}, \widetilde{C} \in X$
        $\widetilde{d}(\widetilde{A},\widetilde{B})= \bigl(\underline{\lambda}_{AB},d(A,B),\overline{\lambda}_{AB}\bigr)$\\
       $\widetilde{d}(\widetilde{A},\widetilde{C})= \bigl(\underline{\lambda}_{AC},d(A,C),\overline{\lambda}_{AC}\bigr)$\\
       $\widetilde{d}(\widetilde{C},\widetilde{B})= \bigl(\underline{\lambda}_{CB},d(C,B),\overline{\lambda}_{CB}\bigr)$
    Here, $d(A,B)$, $d(A,C)$ and $d(C,B)$ are the euclidean distances of the cores of the fuzzy points $\widetilde{A}$, $\widetilde{B}$ and $\widetilde{C}$ respectively. So we have $d(A,B) \leq d(A,C) + d(C,B)$.
\end{enumerate}
\begin{align*}      
    \widetilde{d}(\widetilde{A},\widetilde{C}) + \widetilde{d}(\widetilde{C},\widetilde{B}) &= \bigl(\underline{\lambda}_{AC},d(A,C),\overline{\lambda}_{AC}\bigr) + \bigl(\underline{\lambda}_{CB},d(C,B),\overline{\lambda}_{CB}\bigr)\\
    &= \bigl(\underline{\lambda}_{AC}+\underline{\lambda}_{CB},d(A,C)+d(C,B),~\overline{\lambda}_{AC}+\overline{\lambda}_{CB}\bigr)\\
    &\succeq \bigl(\underline{\lambda}_{AB},d(A,B),\overline{\lambda}_{AB}\bigr)\\
    &= \widetilde{d}(\widetilde{A},\widetilde{B})
\end{align*}   
Therefore, $ \widetilde{M_d}(\widetilde{A},\widetilde{B},t)$ is also fuzzy metric and consequently,  $(X,\widetilde{M_d},*)$ is a metric space.
\subsection{Fuzzy Hausdorff Distance}
The concept of Hausdorff distance can also be extended to the aforementioned fuzzy metric space. In general, Hausdorff distance is a metric that quantifies the similarity or dissimilarity between two sets considering shape and structure of sets. The \emph{Hausdorff distance} between $A$ and $B$, where, $A, B\subset X$ be nonempty sets, is defined as:  
\begin{equation}
\label{Hausdorff}
    d_H(A,B)\;=\;\max\Big\{\sup_{a\in A}\inf_{b\in B} d(a,b),\;\sup_{b\in B}\inf_{a\in A} d(b,a)\Big\}.
\end{equation}

$d_H(A, B)$ represents the greatest of all distances from a point in one set to the closest point in the other set. Essentially, it quantifies the distance between the two sets.

Now we are trying to define the fuzzy Hausdorff distance, which is an extension of the equation (\ref{Hausdorff}). Fuzzy Hausdorff distance can also be interpreted as the distance between two fuzzy sets. However, in our opinion, the distance between two fuzzy geometrical objects cannot be  as precise as given above.  

Using the above concept together with Definition \ref{fuzzydist_set_pt}, the idea is extended to the fuzzy Hausdorff distance between two fuzzy sets. The first step is to determine the Hausdorff distance  between the cores of fuzzy sets $\widetilde{S}_1$ and $\widetilde{S}_2$. If the cores are $p \in S_1$ and $q \in S_2$ then distance is $d(p,q)=d_H\bigl(S_1^{1},S_2^{1}\bigr)$. Considering two fuzzy numbers  $\widetilde{T}(p)=\bigl(T(p)-\gamma_1,T(p),T(p)+\gamma_2\bigr)$ and $\widetilde{T}(q)=\bigl(T(q)-\gamma_3,T(q),T(q)+\gamma_4\bigr)$ with spreads $\gamma_1,\gamma_2$ and $\gamma_3,\gamma_4$ along the line segment joining $p$ and $q$ respectively. This leads to the corresponding expression for the fuzzy Hausdorff distance $\widetilde{d}_H(S_1,S_2)=\widetilde{d}(\widetilde{T}(p),\widetilde{T}(q))$ using \cite{manriquez2025approach}.

With the help of Proposition \ref{cor1}, the fuzzy Hausdorff distance can also be re-defined. 
Let $\widetilde{A}$ and $\widetilde{B}$ be two fuzzy points with circular bases, and $\widetilde{A}^{1}=(a_1,a_2)$ and $\widetilde{B}^{1}=(b_1,b_2)$ are the cores then the Hausdorff distance at $\alpha=1$ is $d_H\bigl(\widetilde{A}^1,\widetilde{B}^1\bigr)$, which is a crisp distance. So,  a line $l_\psi$ containing $(a_1,a_2)$ and $(b_1,b_2)$ exists, that makes an angle $\psi$ with the $x$-axis [see figure \ref{fig:theta_psi}]. Let $l_\theta$ be a line segment containing $(a_1,a_2)$ then $l_\theta$ contains $(b_1,b_2)$ when $\theta~=\psi$ [see Proposition \ref{cor1}]. Thus two inverse points, $\underline{A}_\psi^\alpha$, $\overline{B}_\psi^\alpha$ and $\overline{A}_\psi^\alpha$, $\underline{B}_\psi^\alpha$ with respect to the $\widetilde{A}$ and $\widetilde{B}$ at each $\alpha$-cuts for the angle $\psi$ can be found. Now let $\widetilde{a} \in \widetilde{A}$ and $\widetilde{b} \in \widetilde{B}$ be two continuous fuzzy numbers along the line $l_\psi$, then $\widetilde{a}^\alpha=\bigl(T(\underline{A}_\psi^\alpha),T(a_1,a_2),T(\overline{A}_\psi^\alpha)\bigr)$ and $\widetilde{b}^\alpha=\bigl(T(\underline{B}_\psi^\alpha),T(b_1,b_2),T(\overline{B}_\psi^\alpha)\bigr)$, using the mapping $T$ Eq.\ref{transformation}. For each $\alpha$,  $\widetilde{d}_H^\alpha\bigl(\widetilde{A}^{\alpha}, \widetilde{B}^{\alpha}\bigr)= \widetilde{d}^\alpha\bigl(\widetilde{a}^{\alpha}, \widetilde{b}^{\alpha}\bigr)$ and $\bigcup_{\alpha=0}^{1}  \widetilde{d}_H^{\alpha}\bigl(\widetilde{A}^{\alpha},\widetilde{B}^{\alpha}\bigl)=\bigcup_{\alpha=0}^{1}  \widetilde{d}^{\alpha}\bigl(\widetilde{a}^{\alpha},\widetilde{b}^{\alpha}\bigl)$, that is nothing but the fuzzy distance between two fuzzy points as well. Then the fuzzy Hausdorff distance, denoted by $\widetilde{d}_H(\widetilde{A}, \widetilde{B})$, between two fuzzy points $\widetilde{A}$ and $\widetilde{B}$ is $\widetilde{d}_H(\widetilde{A}, \widetilde{B})=\bigcup_{\alpha=0}^{1}  \widetilde{d}^{\alpha}\bigl(\widetilde{a}^{\alpha},\widetilde{b}^{\alpha}\bigl)~= \widetilde{d}(\widetilde{a}, \widetilde{b})$. Mathematically, the definition is as follows:

\begin{definition}[Fuzzy Hausdorff distance]
The fuzzy Hausdorff distance $\widetilde{d}_H(\widetilde{A},\widetilde{B})$ between two fuzzy points $\widetilde{A}$ and $\widetilde{B}$ in $(X,\widetilde{M_d},*)$ may be defined with its membership function
\begin{center}    $(\widetilde{d}_H^\alpha~|~\widetilde{d}_H)=sup\{\alpha: d=d(u,v),
     ~u \in \widetilde{a}^0, ~v \in \widetilde{b}^0 $ \text{ are inverse points of the continuous fuzzy numbers} $\widetilde{a}$ \text{ and } $\widetilde{b}$ \\ \text{ along the line joining } $\widetilde{A}^1$  \text{and}  $\widetilde{B}^1$  \text{ such that }  $\mu(u\mid \widetilde{a})=\mu(v\mid \widetilde{b})=\alpha$ \text{ and } $d\bigl(\widetilde{a}^1,\widetilde{b}^1\bigr)~=d_H\bigl(\widetilde{a}^1,\widetilde{b}^1 \bigr)$\big\}\\
     \text{and}\\ 
     $\widetilde{d}_H(\widetilde{A},\widetilde{B})~= \bigcup_{\alpha=0}^{1}  \widetilde{d}_H^\alpha\bigl(\widetilde{a}^\alpha,\widetilde{b}^\alpha\bigl)$
\end{center}
\end{definition}

\begin{example}
Considering example
    \ref{fuzzy_dist_ex} for fuzzy points \(\widetilde{A}\), \(\widetilde{B}\) and fuzzy Hausdorff distance is to be defined.
    Then take the line segment $l_\psi: \frac{x-1}{\cos{\psi}}=\frac{y}{\sin{\psi}}$ which contains $(1,0)$ and makes an angle $\psi$ with the $x$-axis. It contains $(5,2)$ when $ \psi~= \tan^{-1}\frac{2-0}{5-1}~=0.46365$. And it includes the line segment that defines the Euclidean distance between $\widetilde{A}^1$ and $\widetilde{B}^1$. So, for the angle $\psi$ we have, [see Figure \ref{fig:theta_psi}],
    \[l_\psi : x-2y=1\]
   Now, if $\widetilde{P}_A=(0.118033989,1.118033989,2.118033989)$ be fuzzy number for fuzzy points $\widetilde{A}$ along the line $l_\psi$ and $\widetilde{P}_B=( 4.472135955, 5.59016994,6.708203932)$ be fuzzy number for fuzzy point $\widetilde{B}$ along the line $l_\psi$, then Hausdorff distance between these two fuzzy points is $\widetilde{d}(\widetilde{P}_A,\widetilde{P}_B)$. Then, Fuzzy Hausdorff distance between two fuzzy points \(\widetilde{A}\), \(\widetilde{B}\) is

$\widetilde{d}_H(\widetilde{A}, \widetilde{B}) = (2.35379606, 4.47213595, 6.59016994)$,\\
which is same as the fuzzy distance between two fuzzy points using Definition \ref{def_fuzzy_distance} i.e.\ $\widetilde{d}_H(\widetilde{A}, \widetilde{B}) \approx \widetilde{d}(\widetilde{A},\widetilde{B})$.
\end{example}


\section{Fuzzy equidistant set}\label{sec_Redefine_proposal}
In this section, we introduce the central concept of this paper: the \emph{fuzzy equidistant set}. These foundational definitions provide the analytical scaffold for the topological and stability analyses presented in the following sections.

In general, for two nonempty sets $A, B\subset X$, mid set \cite{wilkerEquidistantSetsTheir1975} determined by $A$ and $B$ is defined to be
\begin{center}
 $E(A,B)=\left\{x\in X : d(x,A)=d(x,B)\right\}$   
\end{center}
 where  $d(x,A)=\inf\left\{d(x,a):\;a\in A\right\}$. The sets $A$ and $B$ are called focal sets of $E(A,B)$.

The set $E(A,B)$ can be defined in fuzzy metric space. A preliminary work can be found in \cite{manriquez2025approach}. If the given two sets are in fuzzy metric space, since it is derived from a fuzzy metric, fuzzy equidistant set will be fuzzy rather than a crisp subset of $X$ consisting of all points equidistant from the focal sets $A$ and $B$. A natural fuzzification of this construction may be provided by associating a membership function that gradually decreases as points deviate from the locus of exact equidistance, thereby providing a genuine fuzzy representation of the equidistant set.

In this section, we redefine the definition of fuzzy midset in closed form. Fuzzy midset can also be viewed as fuzzy equidistant set. Fuzzy equidistant set can be viewed in two ways. Firstly, considering all $\widetilde{x} \in X$, which is a fuzzy point, such that the fuzzy distance between $\widetilde{x}$ and $\widetilde{A}$ is equal to fuzzy distance between $\widetilde{x}$ and $\widetilde{B}$ using Definition \ref{def_fuzzy_distance}. This is the standard way to view the fuzzy equidistant set. Otherwise, we can also find $\widetilde{x} \in X$ such that the fuzzy degree of closeness between $\widetilde{x}$ and $\widetilde{A}$ is equal to the fuzzy degree of closeness between $\widetilde{x}$ and $\widetilde{B}$. Here we find the equidistant set using equation \ref{eq:Re_def_fuzzy_metric} and definition \ref{def_fuzzy_distance} respectively. In this section it has been proved that that both techniques give the same result.

\begin{definition}[Redefined fuzzy  midset]
Let $\widetilde{A}$ and $\widetilde{B}$ be a fuzzy sets on $\mathbb{R}^2$. The fuzzy midset between $\widetilde{A}$ and $\widetilde{B}$, denoted by $\widetilde{\mathcal{M}}(\widetilde{A},\widetilde{B})$ may be redefined by the collection of fuzzy points $\widetilde{P}$, such that $\widetilde{d}(\widetilde{P},\widetilde{A})=\widetilde{d}(\widetilde{P},\widetilde{B})$ and $\widetilde{P}^{\alpha} \in \mathcal{M}\bigl(\widetilde{A}^{\alpha},\widetilde{B}^{\alpha}\bigr) ~~\forall\alpha \in [0,1] $.
\end{definition}

To find the fuzzy midset or say fuzzy equidistant set, we find the midset of each $\alpha$-cut of the fuzzy focal sets. For the easy of notation we say $\widetilde{E}^{\alpha} = \mathcal{M}\bigl(\widetilde{A}^{\alpha},\widetilde{B}^{\alpha}\bigr)$. So  $\widetilde{E}(\widetilde{A},\widetilde{B}) = \bigcup_{\alpha=0}^{1} \mathcal{M}\bigl(\widetilde{A}^{\alpha},\widetilde{B}^{\alpha}\bigr)$. So for two nonempty sets $\widetilde{A}$ and $\widetilde{B}$ of a fuzzy metric space $(X,\widetilde{d},L,R)$, the fuzzy equidistant set determined by $\widetilde{A}$ and $\widetilde{B}$ is defined to be 
\begin{equation}
\begin{aligned}
\widetilde{E}(\widetilde{A},\widetilde{B})&=\left\{\widetilde{x}\in X : \widetilde{d}(\widetilde{x},\widetilde{A})=\widetilde{d}(\widetilde{x},\widetilde{B})\right\}\\
&= \left\{\bigcup_{\alpha=0}^{1}\widetilde{x}^{\alpha}\in X : \widetilde{d}^{\alpha}\bigl(\widetilde{x}^{\alpha},\widetilde{A}^{\alpha}\bigl)~=~\widetilde{d}^{\alpha}\bigl(\widetilde{x}^{\alpha},\widetilde{B}^{\alpha}\bigl)\right\}\\
&= \left\{\bigcup_{\alpha=0}^{1}\widetilde{x}^{\alpha}\in X : t + \widetilde{d}^{\alpha}\bigl(\widetilde{x}^{\alpha},\widetilde{A}^{\alpha}\bigl)~=~ t + \widetilde{d}^{\alpha}\bigl(\widetilde{x}^{\alpha},\widetilde{B}^{\alpha}\bigl)\right\}\\
&= \left\{\bigcup_{\alpha=0}^{1}\widetilde{x}^{\alpha}\in X : \frac{t}{t + \widetilde{d}^{\alpha}\bigl(\widetilde{x}^{\alpha},\widetilde{A}^{\alpha}\bigl)}=\frac{t}{t + \widetilde{d}^{\alpha}\bigl(\widetilde{x}^{\alpha},\widetilde{B}^{\alpha}\bigl)}\right\}\\
&= \left\{\widetilde{x}\in X : \frac{t}{t + \bigcup_{\alpha=0}^{1}  \widetilde{d}^{\alpha}\bigl(\widetilde{x}^{\alpha},\widetilde{A}^{\alpha}\bigl)}=\frac{t}{t + \bigcup_{\alpha=0}^{1}  \widetilde{d}^{\alpha}\bigl(\widetilde{x}^{\alpha},\widetilde{B}^{\alpha}\bigl)}\right\}\\
&= \left\{\widetilde{x}\in X : \widetilde{M}_d(\widetilde{x},\widetilde{A},t)=\widetilde{M}_d(\widetilde{x},\widetilde{B},t)\right\}
\label{fuzzy_equi_dist_set}
\end{aligned}
\end{equation}

The sets $\widetilde{A}$ and $\widetilde{B}$ are called fuzzy focal sets of $\widetilde{E}(\widetilde{A},\widetilde{B})$. The fuzzy equidistant set $\widetilde{E}(\widetilde{A},\widetilde{B})$ can be defined by its membership function:
\begin{equation}  
\begin{aligned}
\mu_{\widetilde{E}}(x) &= \Biggl\{\,\bigcup_{\alpha=0}^{1}\alpha : \widetilde{d}_\alpha \bigl(
\widetilde{x}^{\alpha}
,\widetilde{A}^{\alpha}\bigr)=\widetilde{d}_\alpha \bigl(
\widetilde{x}^{\alpha}
,\widetilde{B}^{\alpha}\bigr);\quad \widetilde{x} \in X\,\Biggr\}\\
&= \Biggl\{\,\bigcup_{\alpha=0}^{1}\alpha : \widetilde{M}\bigl(
\widetilde{x}^{\alpha}
,\widetilde{A}^{\alpha},t\bigr)=\widetilde{M}\bigl(
\widetilde{x}^{\alpha}
,\widetilde{B}^{\alpha},t\bigr);\quad \widetilde{x} \in X\ ~\&~~ t>0\Biggr\}
\end{aligned}
\end{equation}

where $\widetilde{M}(\widetilde{x},\widetilde{A},t)$ and $\widetilde{M}(\widetilde{x},\widetilde{B},t)$ denote the fuzzy proximities to the fuzzy focal sets $\widetilde{A}$ and $\widetilde{B}$, respectively. Note that the function $\mu_{\widetilde{E}}(x)$ is continuous. Thus, for each $\alpha \in [0,1]$, the $\alpha$-cut is given by
\begin{equation}
\widetilde{E}(\alpha) := \{x \in X : \mu_{\widetilde{E}}(x) \geq \alpha\} 
 = \mu_{\widetilde{E}}^{-1}([\alpha,1]).
\end{equation}.
Thus, if the metric space $(X,\widetilde{M}_d,*)$ is considered, then a similar result for $\widetilde{E}(\widetilde{A},\widetilde{B})$ is obtained.

\subsection{Invariant of fuzzy equidistant set}

The notion of a fuzzy metric space has been introduced in two different forms: one by Kramosil and Mich\'{a}lek \cite{kramosilFuzzyMetricsStatistical1975} and another by Kaleva and Seikkala \cite{kalevaFuzzyMetricSpaces1984}. In both formulations, however, the underlying set is not specified with sufficient clarity. In each case, the fuzzy distance is defined between two crisp points or sets, with the value depending on $t>0$.

In our work, the fuzzy metric space [see definition \ref{Redefined Fuzzy metric space}] is redefined in closed form considering fuzzy distance as fuzzy set, so the fuzzy distance between two fuzzy points [see definition \ref{def_fuzzy_distance}] is itself a fuzzy metric. In this section, it will be shown that the fuzzy equidistant set obtained using both the definitions of fuzzy metric spaces is invariant with the proposed and redefined fuzzy metric space. For convenience, and to distinguish between the two approaches to finding the equidistant set, we denote the sets obtained by the first and second approaches as $\widetilde{E}_{\widetilde{M}}(\widetilde{A},\widetilde{B})$ and $\widetilde{E}_{\widetilde{d}}(\widetilde{A},\widetilde{B})$, respectively.

If the fuzzy metric space  $\left(X,\widetilde{M}_d,*\right)$ is of first type, then for given nonempty subsets $\widetilde{A}, \widetilde{B} \subset X$, the fuzzy equidistant set determined by $\widetilde{A}$ and $\widetilde{B}$ is defined by
\begin{equation}
\widetilde{E}_{\widetilde{M}}(\widetilde{A},\widetilde{B})=\left\{\widetilde{x}\in X : \widetilde{M}_d(\widetilde{x},\widetilde{A},t)=\widetilde{M}_d(\widetilde{x},\widetilde{B},t)\right\}.
\label{ref:fes_KM}
\end{equation}
where the proximity of $x$ to the set $A$ is considered as
\begin{equation*}
\widetilde{M}_d(\widetilde{x},\widetilde{A},t) :=  \frac{t}{t + \widetilde{d}(\widetilde{x},\widetilde{A})}.
\end{equation*}

Again for the second type fuzzy metric space  $\left(X,\widetilde{d},L,R\right)$, if given nonempty subsets are $\widetilde{A}, \widetilde{B} \subset X$, the fuzzy equidistant set is determined by $\widetilde{A}$ and $\widetilde{B}$,
\begin{equation}
\widetilde{E}_{\widetilde{d}}(\widetilde{A},\widetilde{B})=\left\{\widetilde{x}\in X : \widetilde{d}(\widetilde{x},\widetilde{A})=\widetilde{d}(\widetilde{x},\widetilde{B})\right\}.
\label{ref:fes_KS}
\end{equation}

Therefore, from Eq.\ref{ref:fes_KM} and Eq.\ref{ref:fes_KS}, 
\begin{align*}
\widetilde{E}_{\widetilde{M}}(\widetilde{A},\widetilde{B}) &= \left\{\widetilde{x} \in X : \frac{t}{t + \widetilde{d}(\widetilde{x},\widetilde{A})} = \frac{t}{t + \widetilde{d}(\widetilde{x},\widetilde{B})}\right\}\\
&= \left\{\widetilde{x} \in X : \frac{1}{t + \widetilde{d}(\widetilde{x},\widetilde{A})} = \frac{1}{t + \widetilde{d}(\widetilde{x},\widetilde{B})}\right\}\\
&= \left\{\widetilde{x} \in X : t+\widetilde{d}(\widetilde{x},\widetilde{A}) = t+\widetilde{d}(\widetilde{x},\widetilde{B})\right\}\\
&= \left\{\widetilde{x} \in X : \widetilde{d}(\widetilde{x},\widetilde{A}) = \widetilde{d}(\widetilde{x},\widetilde{B})\right\}\\
&= \widetilde{E}_{\widetilde{d}}(\widetilde{A},\widetilde{B}).
\end{align*}
from Eq.\ref{ref:fes_KS},
\begin{align*}
\widetilde{E}_{\widetilde{d}}(\widetilde{A},\widetilde{B}) 
&= \left\{\widetilde{x} \in X : \widetilde{d}(\widetilde{x},\widetilde{A}) = \widetilde{d}(\widetilde{x},\widetilde{B})\right\}\\
&= \widetilde{E}_{\widetilde{M}}(\widetilde{A},\widetilde{B}). \text{  [see Eq.\ref{fuzzy_equi_dist_set}]}
\end{align*}
The same 

\begin{example}
     Let $\widetilde{A}(0,0)$  and $\widetilde{B}(5,0)$ be two fuzzy points with membership functions as 
    \begin{center}
        $\mu_{1}((x,y)|\widetilde{A}(0,0)) = \begin{cases}
        1-\sqrt{(x-0)^{2}+(y-0)^{2} } & if\  (x-0)^{2}+(y-0)^{2} \leq 4 \\
        0 & otherwise \\
    \end{cases}$
    \end{center}

    \begin{center}
        $\mu_{2}((x,y)|\widetilde{B}(5,0)) = \begin{cases}
        1-\sqrt{(x-5)^{2}+(y-0)^{2} } & if\  (x-5)^{2}+(y-0)^{2} \leq 4 \\
        0 & otherwise \\
    \end{cases}$
    \end{center}
    \label{fuzzy_equ_dist_ex1}
\end{example}

Here fuzzy points $\widetilde{A}$ and $\widetilde{B}$ can be expressed as 
\begin{align*}
    &\widetilde{A} : x^2 + y^2 = 4(1-\alpha)^2\\
    &\widetilde{B} : (x-5)^2 + y^2 = 4(1-\alpha)^2
\end{align*}

When $\alpha = 1$, then it indicates the core of the fuzzy points. Which are basically the crisp points in $\mathbb{R}^2$. Then the equidistant line will be the bisector line of the points $(0,0)$ and $(5,0)$. Which is the $x=2.5$ line parallel to the $y$-axis.\\
 
Given $(p,q) \in \mathbb{R}^2$, for distinct values of $\alpha$, the fuzzy equidistant set can be determined by equating the expressions: $
\sqrt{p^2 + q^2} - 2(1-\alpha) = \sqrt{(p-5)^2 + q^2} - 2(1-\alpha)
$. This simplifies to the line $p=2.5$. Consequently, for varying $\alpha$, the equidistant set is defined by $x=2.5$. It is noteworthy that if the spreads of two fuzzy points are circles with identical radii, the resulting equidistant set will be a crisp line. Next example considers the points with different radii.

\begin{figure}[h]
\begin{center} 
\begin{subfigure}{0.4\textwidth}
\includegraphics[width=0.9\linewidth]{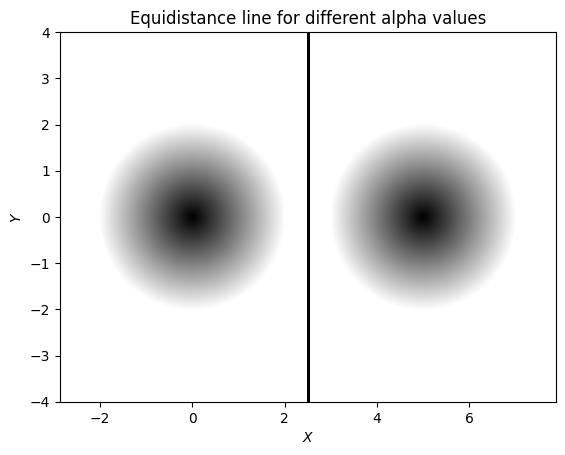} 
\caption{Example \ref{fuzzy_equ_dist_ex1}}
\label{fig:FES1}
\end{subfigure}
\begin{subfigure}{0.4\textwidth}
\includegraphics[width=0.9\linewidth]{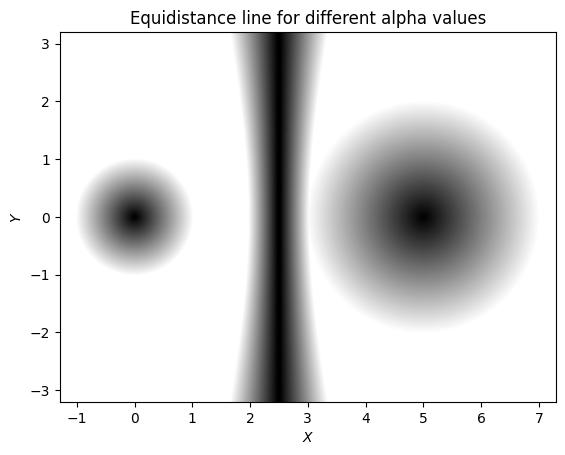}
\caption{Example \ref{fuzzy_equ_dist_ex2}}
\label{fig:FES2}
\end{subfigure}

\caption{In \ref{fig:FES1} the fuzzy points $\widetilde{A}$ and $\widetilde{B}$ are symmetric, so the fuzzy equidistant line is a crisp set. But in \ref{fig:FES2} fuzziness of $\widetilde{B}$ is larger than the spread of $\widetilde{A}$ and they are not overlapping. So it gives fuzzy equidistant set looks like a hyperbola.}
\label{fig:FES}
\end{center}
\end{figure}

\begin{example}
    Let $\widetilde{A}(0,0)$  and $\widetilde{B}(5,0)$ be two fuzzy points with membership functions as 
    \begin{center}
        $\mu_{1}((x,y)|\widetilde{A}(0,0)) = \begin{cases}
        1-\sqrt{(x-0)^{2}+(y-0)^{2} } & if\  (x-0)^{2}+(y-0)^{2} \leq 1 \\
        0 & otherwise \\
    \end{cases}$
    \end{center}

    \begin{center}
        $\mu_{2}((x,y)|\widetilde{B}(5,0)) = \begin{cases}
        1-\sqrt{(x-5)^{2}+(y-0)^{2} } & if\  (x-5)^{2}+(y-0)^{2} \leq 4 \\
        0 & otherwise \\
    \end{cases}$
    \end{center}
    \label{fuzzy_equ_dist_ex2}
\end{example}

Fuzzy points $\widetilde{A}$ and $\widetilde{B}$ can be express as 
\begin{align*}
    &\widetilde{A} : x^2 + y^2 = (1-\alpha)^2\\
    &\widetilde{B} : (x-5)^2 + y^2 = 4(1-\alpha)^2
\end{align*}

Similar to Ex.\ref{fuzzy_equ_dist_ex1} we get $x=2.5$ line for $\alpha=1$. And, for given $(p,q) \in \mathbb{R}^2$ the fuzzy equidistant set can be found from 
\begin{align*}
    \sqrt{p^2 + q^2} - (1-\alpha) &= \sqrt{(p-5)^2 + q^2} - 2(1-\alpha) \implies
\frac{(2p-5)^2}{(1-\alpha)^2} - \frac{4q^2}{25 - (1-\alpha)^2}~=~1
\end{align*}

This gives an equation of a hyperbola.\\
Now let us consider different values of $\alpha$, to generate distinct hyperbolic lines. For $\alpha=1$, we get the cores of the fuzzy points, which gives the $x = 2.5$ line as the equidistant set. Conversely, when $\alpha = 0$, this becomes the equidistant line between two circles, represented by the equation $\frac{(x-5/2)^2}{1/4} - \frac{y^2}{6} = 1$. Thus, for $\alpha \in [0,1]$ different fuzzy equidistant sets are obtained.

Next we give examples of different cases where fuzziness of two fuzzy points may be separated, partially overlapped or fully overlapped. We also try to analyze the fuzzy equdistant set. 

\begin{example}
Let $\widetilde{A}(a_1,a_2)$ and $\widetilde{B}(b_1,b_2)$ be two fuzzy points with circular spreads with different sizes. And let the memberships functions are
\begin{center}
        $\mu_{1}((x,y)|\widetilde{A}(a_1,a_2)) = \begin{cases}
        1-\sqrt{(x-a_1)^{2}+(y-a_2)^{2} } & if\  (x-a_1)^{2}+(y-a_2)^{2} \leq r_1 \\
        0 & otherwise \\
    \end{cases}$
    \end{center}

    \begin{center}
        $\mu_{2}((x,y)|\widetilde{A}(b_1,b_2)) = \begin{cases}
        1-\sqrt{(x-b_1)^{2}+(y-b_2)^{2} } & if\  (x-b_1)^{2}+(y-b_2)^{2} \leq r_2 \\
        0 & otherwise \\
    \end{cases}$
    \end{center}
\end{example}

The nature of the fuzzy equidistant sets with respect to different positions of the fuzzy spreads is examined. The necessary conditions for the different cases are provided. \\

Considering, $d_1=\sqrt{(p-a_1)^{2}+(q-a_2)^{2} }$ and $d_2=\sqrt{(p-b_1)^{2}+(q-b_2)^{2} }$ for $(p,q) \in \mathbb{R}^2$, the different cases are as follows:\\\\
\textbf{Case 1: Fuzzy points are non-overlapping}\\

For each $\alpha \in [0,1]$ the equidistant set is
$d_1 - r_1(1-\alpha) = d_2 - r_2(1-\alpha)$. This equation always gives a hyperbolic typed equidistant set for any distinct $r_1$ and $r_2$. When spreads of the fuzzy points are same, i.e. $r_1 = r_2$ it implies that $d_1 = d_2$. Consequently a precise straight line is obtained, and this serves as the bisector of $\widetilde{A}^1$ and $\widetilde{B}^1$.\\\\
\textbf{Case 2: Fuzzy Points are partially overlapping}\\

For each $\alpha \in [0,1]$ the equidistant set is either $d_1 - r_1(1-\alpha) = d_2 - r_2(1-\alpha)$ or, $d_1 - r_1(1-\alpha) =  r_2(1-\alpha) - d_2$. The first set acts like case-1 that is hyperbolic typed equidistant set. For the second set, when $\alpha=0$ the fuzzy points are partially overlapping and resulting an elliptic line. As the value of $\alpha$ approaches $1$, the fuzzy spreads are reducing and after a certain value of $\alpha=n~(say)$ spreads cease to overlap. Subsequently, from $\alpha \in [n,1]$ it behaves like case-1. So, for $\alpha \in [0,n]$ both hyperbolic and elliptical equidistant set is formed and for $\alpha \in [n,1]$ hyperbolic typed equidistant set is generated.\\\\
\textbf{Case 3: Fuzzy Points are fully overlapping}\\

For each $\alpha \in [0,n_1]$ the equidistant set is
$d_1 - r_1(1-\alpha) =  r_2(1-\alpha) - d_2$, so always it gives elliptical equidistant set. But when the value of $\alpha$ is moving toward $1$ the fuzzy spreads diminish, and for $\alpha=n_1 ~(say)$, the fuzzy points become partially overlapping, then it behaves like case-2. So, it gives two types of equidistant sets when $\alpha \in [n_1,1]$, either $d_1 - r_1(1-\alpha) = d_2 - r_2(1-\alpha)$ or,
$d_1 - r_1(1-\alpha) =  r_2(1-\alpha) - d_2$. Similar to case-2, for the values of $\alpha \in [n_1,n_2]$ (say), both hyperbolic-type and elliptical equidistant sets are generated. While, hyperbolic-type equidistant set for the value of $\alpha \in [n_2,1]$ is generated. If two fuzzy points have the same core but different spreads it always gives elliptical equidistant set.\\

The general equation of the equidistant set is $Ax^2 + 2Hxy + By^2 + 2Gx + 2Fy + C = 0$.
Now according to the relative positions of the fuzzy points and fuzzy spreads above conic sections as a equidistant set are achieved.

If $ABC + 2FGH - AF^2 - BG^2 - CH^2=0$ and $AB-H^2=0$ the equidistant set is a straight line. This condition holds only if circular spreads are equal, resulting in a distinct and precise straight line.
Conversely, if $ABC + 2FGH - AF^2 - BG^2 - CH^2 \neq 0$, the equidistant set is either an elliptic or a hyperbolic. If $AB - H^2 > 0$, the equidistant set is an elliptic; otherwise, it is hyperbolic. But $d_1 - r_1(1-\alpha) = d_2 - r_2(1-\alpha)$ is obtained by computing the inverse points and $d_1 - r_1(1-\alpha) = r_2(1-\alpha) - d_2$ is obtained by computing the same points. Figures of the different cases for better understanding are provided later. However, by the definition of fuzzy distance [see Definitions \ref{def_fuzzy_distance} and \ref{fuzzydist_set_pt}] the hyperbolic type equidistant set is accepted.

\begin{figure}[h!]
\centering{\includegraphics[scale=0.7]{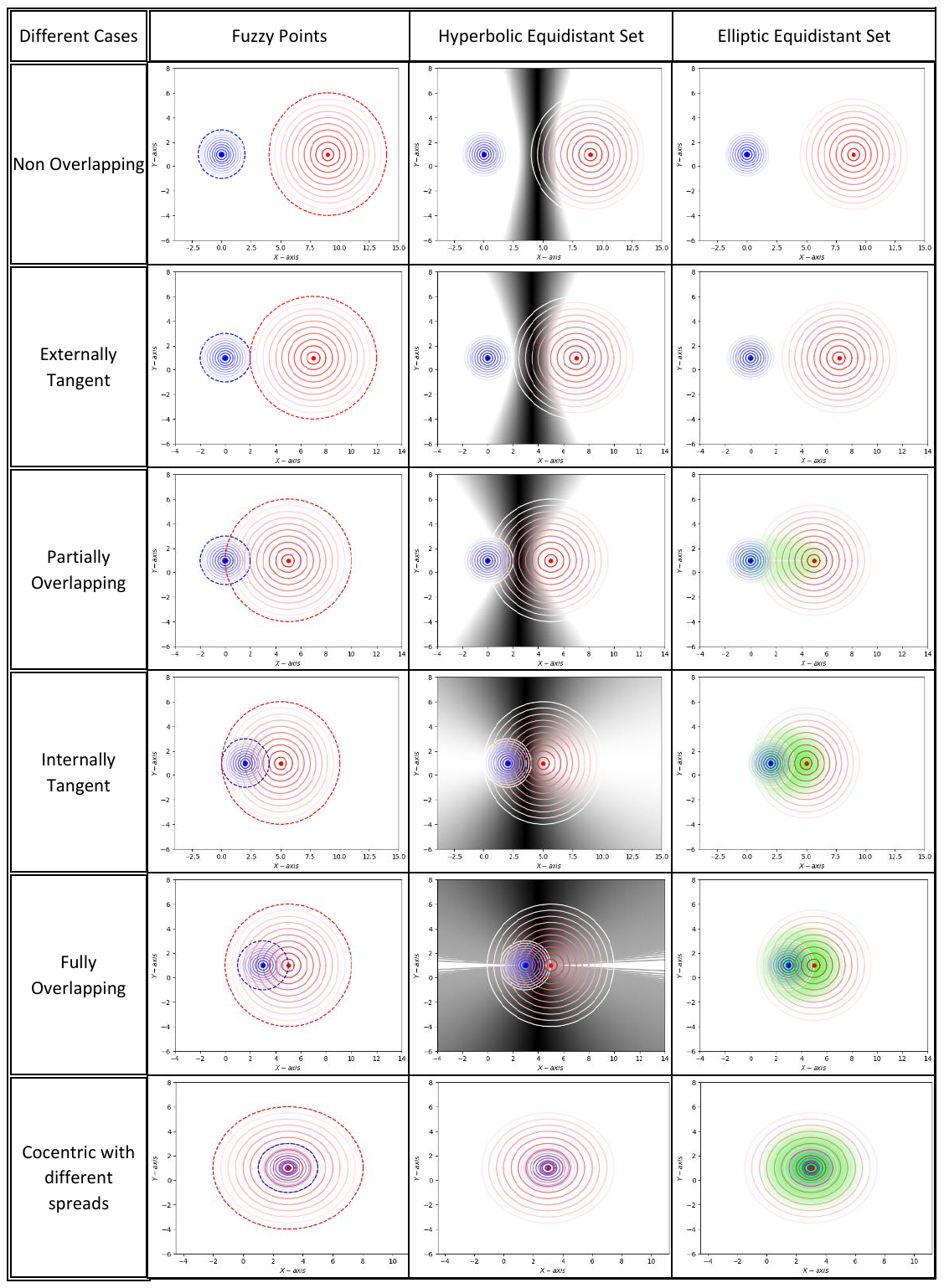}}
	\caption{Fuzzy equidistant set between different fuzzy points $\widetilde{A}$ and $\widetilde{B}$}
	\label{FPs_table} 
\end{figure}

\begin{landscape}
\begin{center}    
\begin{tabular}{ |p{2cm}||p{3.2cm}|p{1.7cm}|p{14cm}|} 
 \hline
 \multicolumn{4}{|c|}{ {\textbf{Nature of the Equidistant Set Between Two fuzzy points with circular spreads of different sizes}}}\\
 \hline
{ \makecell{\textbf{Relative}\\\textbf{Positions}}}& {\makecell{\textbf{Condition}}} &{ \makecell{\textbf{Nature}}}&{ \makecell{\textbf{Remarks}}}\\
 \hline
 \hline
{\makecell{ fuzzy spreads\\are non\\overlapping}}   & { \makecell{$d_c > (r_1+r_2)(1-\alpha)$,\\ $\forall \alpha \in [0,1]$}}    & { \makecell{hyperbolic}} &   { For all $\alpha \in [0,1]$ fuzzy spreads are separated and generates a fuzzy hyperbolic equidistant set.}\\
 \hline
   { \makecell{fuzzy spreads\\are externally\\tangent}}&   { \makecell{$d_c \geq (r_1+r_2)(1-\alpha)$,\\ $\forall \alpha \in [0,1]$}}  & { \makecell{hyperbolic}}   &{ for $\alpha=0$ one circle is externally tangent to another circle. And $\forall \alpha \in (0,1]$, fuzzy spreads are separated. For all cases, the equidistance set is hyperbolic.}\\
 \hline
  { \makecell{fuzzy spreads\\ are partially\\overlapping}}  & { \makecell{$(r_1-r_2)(1-\alpha)<d_c$ \\ $< (r_1+r_2)(1-\alpha)$, \\$\forall \alpha \in [0,n]$}}    &  { \makecell{elliptic\\+\\hyperbolic}} &  { for $\alpha \in [0,n]$ the spreads of fuzzy points are partially overlapping, so an elliptical equidistant set will generate $\forall \alpha \in [0,n]$ for the same points. And, a hyperbolic-type equidistant set will generate $\forall \alpha \in [0,n]$ for the inverse points.} \\
 &
 { \makecell{
 $d_c \geq (r_1+r_2)(1-\alpha)$,\\ $\forall \alpha \in [n,1]$}} & { \makecell{hyperbolic}}& { Here $n$ is any number between 0 and 1 and depends on the fuzzy spreads. As $\forall \alpha \in [n,1]$ the spreads start to separate and make a hyperbolic equidistant set.}\\
 \hline
  { \makecell{fuzzy\\spreads are\\internally\\tangent}} & { \makecell{$(r_1-r_2)(1-\alpha)\leq d_c$ \\ $< (r_1+r_2)(1-\alpha)$, \\$\forall \alpha \in [0,n]$}}    &  { \makecell{elliptic\\+\\hyperbolic}} &  { for $\alpha \in [0,n]$ the spreads of fuzzy points are internally tangent or partially overlapping and $\forall \alpha \in [0,n]$ both elliptical and hyperbolic-type equidistant set will be formed for same points and inverse points respectively.}\\
 &
  { \makecell{
 $d_c \geq (r_1+r_2)(1-\alpha)$,\\ $\forall \alpha \in [n,1]$}} & { \makecell{hyperbolic}}& { Here $n$ is any number between 0 and 1 and depends on the fuzzy spreads. As $\forall \alpha \in [n,1]$ the spreads start to separate and make a fuzzy hyperbolic equidistant set.}\\
 \hline
     & { \makecell{$d_c \leq (r_1-r_2)(1-\alpha)$,\\ $\forall \alpha \in [0,n_1]$}}    &  { \makecell{elliptic}} &  { For $\alpha \in [0,n_1]$ fuzzy spreads are fully overlapping. And for $\alpha \in [n_1,n_2]$, fuzzy spreads are partially overlapping.}\\
       { \makecell{fuzzy spreads\\are fully\\overlapping}}  & { \makecell{$(r_1-r_2)(1-\alpha)\leq d_c$ \\ $< (r_1+r_2)(1-\alpha)$, \\$\forall \alpha \in [n_1,n_2]$}}   & { \makecell{elliptic\\+\\hyperbolic}} &  { For any numbers between $0$ and $1$, $0 \geq n_1\geq n_2 \geq 1$ fuzzy elliptical equidistant set is obtained for $\alpha \in [0,n_1]$. And for $\alpha \in [n_1,n_2]$, both elliptic and hyperbolic-type equidistant sets are generated for the same points and inverse points respectively.}\\
 
 &
 { \makecell{
 $d_c \geq (r_1+r_2)(1-\alpha)$,\\ $\forall \alpha \in [n_2,1]$}} &  { \makecell{hyperbolic}}& { For $\alpha \in [n_2,1]$ fuzzy spreads start to separate and form fuzzy hyperbolic equidistant set}.\\
 \hline
 { \makecell{Concentric\\fuzzy points}}& { \makecell{$d_c =0, ~~\forall\alpha \in [0,1]$}}  & { \makecell{elliptic}}&{ For all $\alpha \in [0,1]$ fuzzy spreads are co-centric and generates a fuzzy elliptic equidistant set.}\\
 \hline
\end{tabular}
\end{center}
\end{landscape}\

\section{Conclusions and comments}\label{sec_conclusion}

This paper formally introduces the notion of a fuzzy equidistant set generated by two non-empty fuzzy subsets of a fuzzy metric space $(X,\widetilde{M},\ast)$.  For a fixed parameter $t>0$, the degree of closeness between two sets are explained as a fuzzy set. We also analyzed that there is less uncertainty in the degree of closeness when t is nearer to 0 or very high. 
Moreover, the fuzzy equidistant set $\widetilde{E}(\widetilde{A},\widetilde{B})$ has been shown to be invariant under fuzzy metric spaces $(X,\widetilde{M},\ast)$ and $(X,\widetilde{d},L,R)$. In the majority of the papers, fuzzy sets were not regarded as subsets of fuzzy metric spaces. Instead, a probabilistic approach was adopted to calculate the distance between two crisp points in fuzzy metric space with respect to a parameter $t>0$.   In this paper, the uncertainty in the degree of closeness is attributed due to fuzziness not randomness, as the distance metric is itself fuzzy. It is worthy to mention that here geometry has been applied for fuzzy convex sets. Similar study may be adopted for non-convex fuzzy sets \cite{chakraborty2024introduction} and for other geometrical structures also.
A systematic study of fuzzy equidistant sets and their topological or metric stability has remained underexplored. A comprehensive study on topology is planned for future work.

\bibliographystyle{plain}
\bibliography{Fuzzy_midset_Fmetric.bib} 

\end{document}